\documentclass{article}
\usepackage{tikz}
\usepackage{geometry}
\usepackage{fancyhdr}
\usepackage{amsmath, amssymb}
\usepackage{graphicx}
\usepackage{color}
\usepackage{verbatim}
\definecolor{bg}{rgb}{0.89, 0.95, 0.71} 
\definecolor{ac}{rgb}{0.50, 0.18, 0.41} 
\definecolor{rc}{rgb}{0.77, 0.57, 0.71} 
\usepackage{listings}
\usetikzlibrary{arrows,intersections}

\lstdefinestyle{latex}{language=TeX,
                       backgroundcolor=\color{bg},
                       basicstyle=\small\ttfamily,
                       frame=leftline,
                       xleftmargin=1.4em,
                       framexleftmargin=.8em}
\lstdefinestyle{cmdline}{
                         }
\usepackage{url}
\usepackage[pdftex,
            bookmarks,
            colorlinks=true,
            linkcolor=black,
            plainpages = false,
            pdfpagemode = UseNone,
            pdfstartview = FitH,
            citecolor = ac, urlcolor = ac, filecolor = ac]{hyperref}

\usepackage[para]{footmisc}
\makeatletter
\long\def\@makefntext#1{\leavevmode
\@makefnmark\nobreak
\hskip.05em\relax#1%
}
\makeatother

\makeatletter 
\def\section{\@startsection {section}{1}{\z@}{2.5ex plus .6ex minus
    .2ex}{1.0ex plus .15ex}{\hspace*{-3em}\Large\bf\color{ac}}}
\def\subsection{\@startsection{subsection}{2}{\z@}{1.5ex plus .3ex minus
   .1ex}{.2ex plus .1ex}{\hspace*{-3em}\bf\large\color{ac}}}
\makeatother
\newlength{\rulelength}
\title{Getting something out of \LaTeX{}}
\author{Jim Hef{}feron}
\usepackage{wrapfig}
\usepackage{subfig}

\RequirePackage{lineno}
\modulolinenumbers[1]

\expandafter\def\expandafter\normalsize\expandafter{%
    \normalsize
    \setlength\abovedisplayskip{5pt}
    \setlength\belowdisplayskip{5pt}
    \setlength\abovedisplayshortskip{-5pt}
    \setlength\belowdisplayshortskip{10pt}
}

\begin{document}

\thispagestyle{empty}

\begin{center}



\vspace*{-1.5cm}

{\large\bf An Exploratory Study on Machine Learning to Couple 

Numerical Solutions of Partial Differential Equations}

\vspace*{0.2cm}

H. S. Tang$^{a,}$\footnote[1]{Correspondence: Hasnong Tang, htang@ccny.cuny.edu}, L. Li$^{a,b}$, M. Grossberg$^c$, Y. J. Liu$^d$, Y. M. Jia$^e$, S. S. Li$^{a,f}$, W. B. Dong $^a$ 

$^a$ Civil Engineering Department, City College of New York 

City University of New York, New York 10031, USA

$^b$ Department of Control Science and Engineering 

University of Shanghai for Science and Technology, Shanghai 200093, China

$^c$ Department of Computer Science, City College of New York 

City University of New York, New York 10031, USA

$^d$ School of Mathematics, Georgia Institute of Technology, Georgia 30332, USA 

$^e$ School of Automation Science and Electrical Engineering 

Beihang university, Beijing 10083, China

$^f$ School of Hydraulic and Hydroelectric Engineering 

Xi’an University of Technology, Shaanxi 710048, China

\end{center}

\vspace*{-0.1cm}

\begin{center}
\parbox{0.9\hsize}
{\small

{\sf

{\bf ABSTRACT} \ 
As further progress in the accurate and efficient computation of coupled partial differential equations (PDEs) becomes increasingly difficult,  it has become highly desired to develop new methods for such computation. In deviation from conventional approaches, this short communication paper explores a computational paradigm that couples numerical solutions of PDEs via machine-learning (ML) based methods, together with a preliminary study on the paradigm. Particularly, it solves PDEs in subdomains as in a conventional approach but develops and trains artificial neural networks (ANN) to couple the PDEs' solutions at their interfaces, leading to solutions to the PDEs in the whole domains. The concepts and algorithms for the ML coupling are discussed using coupled Poisson equations and coupled advection-diffusion equations. Preliminary numerical examples illustrate the feasibility and performance of the ML coupling. Although preliminary, the results of this exploratory study indicate that the ML paradigm is promising and deserves further research. 

\vspace{0.1cm}

{\it Keywords}: Machine learning, artificial neural network, coupled partial differential equations, domain decomposition,  interface zone, transmission condition}
}
\end{center}

\vspace{0.2cm}

\noindent {\bf 1. Introduction}

\vspace{0.1cm}


The research on the coupling of partial differential equations (PDEs) can be traced back to Schwarz's idea of domain decomposition over a hundred years ago \cite{Schwarz1869}, and it attracts more and more attention in the past few decades. For relatively simple PDEs such as advection-diffusion equations, the past findings include slow convergence associate with a regular Dirichlet transmission condition \cite {Califano2018}, a decrease in convergence speed with grid spacing \cite {Gastaldi1989, Gander1998}, and speedup of the convergence by methods such as precondition in iteration matrices and optimization for interface conditions \cite {Martin2004, Gander2007}. For systems of conservation laws, analysis indicates a convergence of numerical solutions to weak solutions in the association of both conservative and non-conservative interface treatments \cite {Berger1987, Tang1999}. In recent years, there is a trend to couple distinct, complex PDEs for simulation of real-world multiscale and multiphysics problems \cite {keyes2012,Tang2020}. Examples are the computation of an acoustic fluid-structure interaction problem by a monolithic balancing domain decomposition method \cite{Minami:2012jp}, simulation of ocean flows by the integration of the Naver-Stokes equations and their hydrostatic versions \cite {Tang2014a, blayo2016, Qu2019a}, and many more. 

Despite the fact that substantial research has been made, progress in theoretical and applied research is limited because of the complexity in coupling PDEs, especially different types of PDEs. The theoretical analysis has mostly been restricted to model/simplified problems, e.g., \cite {Eisenmann2018}, while the applied research deals with realistic problems, but its rigorous foundation is unclear, e.g., \cite {Qu2019a}. Actually, difficulties are encountered in the computation of various problems. For instance, in the computation of hyperbolic systems of conservation laws, numerical oscillations, nonphysical solutions, and multiple solutions may take place \cite {Part1994, Wu1996, Tang1996}. In simulations of ocean flows, especially highly transient flows, problems include artifacts (e.g., numerical oscillations) at the interfaces between PDEs, numerical instability, and slow convergence of the computation  \cite {Tang2014a, Tang2016a, Tang2019}. As a result, our modeling capabilities are impaired and cannot meet the needs of many problems.  

Further progress in coupling PDEs, especially distinct PDEs, via existing approaches becomes increasingly difficult. The difficulties stem from the fact that the coupling is a heterogeneous domain decomposition problem; it involves different PDEs, numerical methods, and computational meshes. The difficulties are hard to overcome in a conventional coupling approach, which directly exchanges PDEs' solutions at their interfaces according to certain interface conditions. Frequently, such interface conditions are complicated, e.g., between distinct PDEs, and, if available, they lack a rigorous foundation in mathematics. As a result, the overall coupled problems' well-posedness is often unclear or even not in presence. When coupling distinct PDEs via a conventional approach, usually application of advanced techniques, e.g., a conservative interface treatment and the Newton-Krylov-Schwarz method \cite {Tang2003,Chen2014b}, are neither straightforward nor effective. 

Now we see the rapid progress of machine learning (ML) in methods and also applications in a wide range of problem domains, including speech recognition, image classification, and superhuman performance on games previously thought intractable such as Go  \cite {Yassin2018, Park2018}. More recently, such progress is taking place also in PDE-based learning of fluid flows \cite {Tompson2017,Raissi2019}. A particular exciting achievement is that, based on status at a current time, deep learning is capable of predicting flow fields at a later time \cite {Lee2017,Wiewel2019}. Analysis indicates the convergence of a deep learning solution to that of the PDE in terms of the numbers of neurons \cite {Sirignan2018}. Now, many powerful, open-source ML tools are available for the public to use, such as Matlab, Pytorch, and Tensorflow \cite {matlabweb,pytorchweb,tensorflowweb}. Interestingly, ML has been utilized to integrate micro-macro fluid flows very recently \cite {Pawar2020}. 

How about coupling PDEs via machine learning? This is a natural idea given the bottleneck difficulties encountered in conventional approaches and inspiration from machine learning's success. If successful, coupling PDEs via ML methods may help us overcome these difficulties and lead to a new avenue to accurate simulation of real-world multiscale multiphysics problems. We have started to explore such an idea since last year, and this short communication paper reports partial investigation and results obtained since then (e.g., \cite {Tang2019b}). Particularly, it presents the concepts and methods of ML-based coupling of Poisson equations as well as of advection-diffusion equations, and it illustrates its feasibility and performance via numerical examples. Although preliminary, it is expected that this short communication could be interesting to the community and will attract its attention to ML coupling of PDEs.

\vspace{0.2cm}

\noindent {\bf 2. Poisson equation} 

\vspace{0.1cm}

Consider the boundary value problem of the Poisson equation as follows: 

\begin {equation}
\label{Poisson}
\left\{\begin{array}{lll}
\Delta u=f, & \  \textbf{x}\in \Omega , \\
u=g, & \  \textbf{x}\in \partial  \Omega ,
\end{array}\right.
\end {equation}

\noindent where $u, f, g\in C^0$. $u $ is the unknown function, and $f $ and $g$ are prescribed functions. 

\begin{wrapfigure}{r}{0.4\textwidth} 
\vspace {-0.5 cm}
\centering
\includegraphics[width=0.99\linewidth]{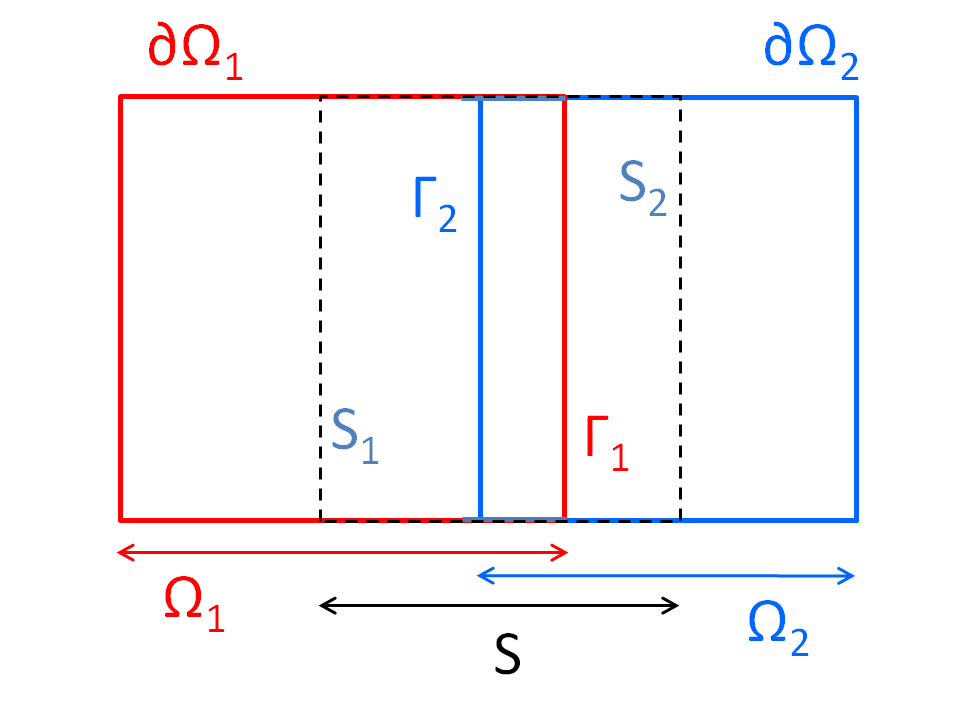} 
\vspace {-0.6 cm}
\caption{\small Domain partition. The dash lines enclose the interface zone, and it contains $\Gamma _1 $ and $\Gamma _2$, the interfaces of the two subdomains.}
\label{domain_division}
\vspace {-0.2 cm}
\end{wrapfigure}

Let domain $\Omega $ be divided into two overlapping subdomains $\Omega _1$ and $\Omega _2$, which are bounded by $\partial \Omega _1$ and $\partial \Omega _2$, respectively, as shown in Fig. \ref {domain_division}. In the figure, $\Gamma _1$ and $\Gamma _2$ are interface of subdomain $\Omega _1 $ and $\Omega _2$, respectively, and they are artificial boundaries for the two subdoamins. An ML Schwarz iteration algorithm to solve problem (\ref {Poisson}) is presented as follows ($m=0,1,2, ...$):

\begin{equation}
\label{ML_Schwarz1}
\left\{\begin{array}{ll}
\Delta u_1^{m+1}=f, & \ \textbf{x}\in \Omega _1, \\
u_1^{m+1}=s^m, & \ \textbf{x} \in \Gamma _1, \\
u_1^{m+1}=g, & \ \textbf{x} \in \partial  \Omega _1/\Gamma _1 ,
\end{array}\right . \ \ \ \ \ 
\left\{\begin{array}{ll}
\Delta u_2^{m+1}=f, & \ \textbf{x}\in \Omega _2, \\
u_2^{m+1}=s^m, & \ \textbf{x} \in \Gamma _2, \\
u_2^{m+1}=g, & \ \textbf{x} \in \partial  \Omega _2/\Gamma _2 ,
\end{array}\right .
\end{equation}

\noindent where $s$ is an ML solution, and it computes the transmission conditions, or interface conditions, on $\Gamma _1 $ and $\Gamma _2$ to facilitate the computation. When the interface conditions are replaced by $u_1^{m+1}=u_2^m$ and $u_2^{m+1}=u_1^{m+1}$, algorithm (\ref {domain_division}) becomes the original alternative Schwarz iteration \cite {Schwarz1869,Toselli2005}.

Let $S$ be an interface zone wrapping the interfaces of the two subdomains and be bounded by $S_1$, $S_2$, and $\partial \Omega _1\cap \partial \Omega _2$ (Fig. \ref {domain_division}). The ML solution, $s$, will be constructed within this zone as a surrogate solution or an approximate solution to the Poisson equations. Let $u$ be the corresponding solution of the Poisson equation, and $u_0 $ be its boundary condition. Then the problem of the ML solution becomes a boundary value problem: Find $\mathcal {M}$ such that

\begin {equation} 
\label{ML_problem1}
\left\{\begin{array}{ll} 
s= \mathcal {M}(u_0)\approx u,  \ \ \ & \textbf{x}\in S \\
s=u_0, \ \ \ & \textbf{x}\in S_1, S_2, \partial \Omega _1\cap \partial \Omega _2
\end{array}\right . 
\end {equation}

\noindent $\mathcal {M}$ may be considered as an operator mimicking the Poisson equation. 

As in a standard approach, the ML solution is constructed as 
\begin {equation} 
\label{interface_solu1}
s=N(t,\textbf {x}, u_0, \boldsymbol {\alpha })
\end {equation}

\noindent where $N$ is a neural network function, and $\boldsymbol {\alpha  }$ denotes its parameters, including weights and biases. In a more general approach, the ML solution may be constructed as a function of $N$ \cite {Lagaris1998,Baymani2015}, for instance, in a numerical example in the following section. Moreover, an ML solution with multi-fidelity will be considered \cite {Kou2019, Liu2019}. 

Since an ML solution, $s$, is a neural network function, $N$, or, a function of $N$, a key to its construction becomes to find the best weights and biases, ${\boldsymbol {\alpha  }}$, within the network. In correspondence to (\ref {interface_solu1}), this leads to a minimization problem: 

\begin {equation} 
\label{loss}
\left\{\begin{array}{ll}
\mbox {Find} \ {\boldsymbol {\alpha  }} \ \mbox {such that} \  Loss(s) \ \mbox{is minimized},\\
\mbox {where} \ Loss(s)=||s-\bar {s}||^2_{2,S}
\end{array}\right. 
\end {equation} 

\noindent in which $Loss$ is the loss function, and $\bar {s} $ is data that approximates the solution of the Poisson equation and trains the network. A main task in finding an ML solution is to efficiently solve the minimization problems, e.g., (\ref {loss}), particularly to search for the best weights and biases in it. Thanks to those researchers who made a breakthrough and developed open-source software packages that can efficiently solve the minimization problem. Available packages include Matlab made by MathWorks \cite {matlabweb}, Tensorflow by the Google Brain Team \cite {tensorflowweb}, and PyTorch primarily by the Facebook's AI Research lab \cite {pytorchweb}. Note that, instead of training via the data, it is possible to train a network via deep learning based on the Poisson equation \cite {Lagaris1998}.

\begin{wrapfigure}{r}{0.4\textwidth} 
\vspace {-0.3 cm}
\centering
\includegraphics[width=0.99\linewidth]{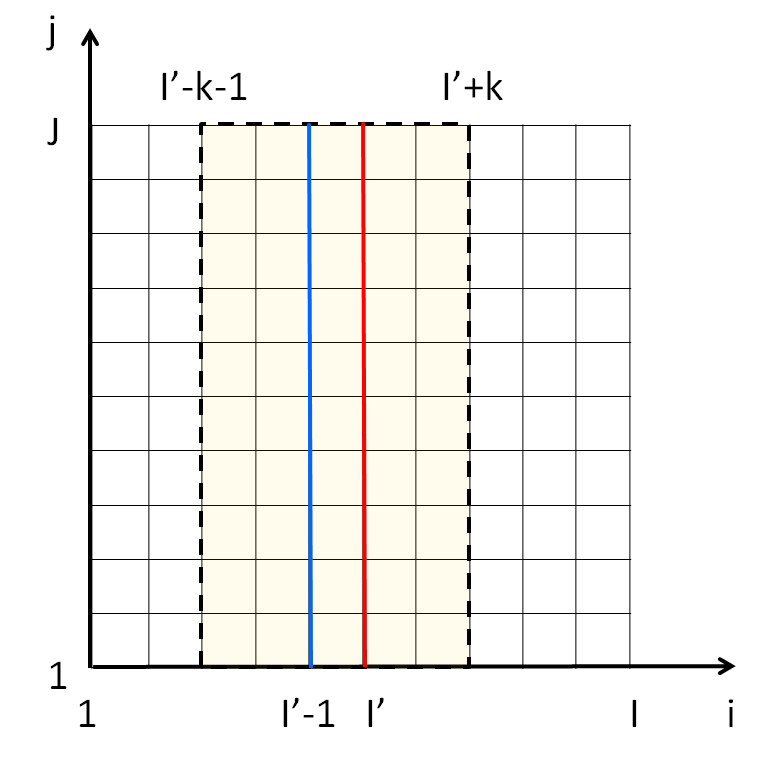} 
\vspace {-0.9 cm}
\caption{\small Schematic representation of the mesh for computation (\ref {ML_numerical}). The red and blue lines are interfaces, and the shadow area enclosed by the dash line is the interface zone.}
\label{mesh}
\vspace {-0.2 cm}
\end{wrapfigure}

In order to illustrate the feasibility and performance of the ML-based coupling methods proposed above, here a numerical example is made.  
Consider a boundary value problem, in which $\Omega : \ 0<x<1,\ 0<y<1$, $f=e^{-x}(x-2+y^3+6y)$, $g(0,y)=y^3, \ g(1,y)=(1+y^3)e^{-1}, \ g(x,0)=xe^{-x}$, and $ g(x,1)=e^{-x}(x+1)$. The problem has an analytical solution \cite {Lagaris1998}: $u_{exact}=e^{-x}(x+y^3)$. Let domain $\Omega $ be discretized by mesh ($x_i,y_j$), $i=1,2,...,I, \ j=1,2, ..., J$. On this grid, let $\Omega _1$ be covered by $i=1,2,...,I', \ j=1,2, ..., J$ and $\Omega _2$ by $i=I'-1,I',...,I, \ j=1,2, ..., J$, as shown in Fig. \ref {mesh}. Moreover, suppose the interface zone is covered by  $i=I_{I'-k-1}, I_{I'-k},..., I_{I'+k}, \ j=1,2, ...,J$, with stencils of $2k+2$ ($k$ is an integer) in the $i$ direction. 

On the grid and by central difference, problem (\ref {ML_Schwarz1}) is discretized as

\begin {equation} 
\label {ML_numerical}
\left\{\begin{array}{ll}
\dfrac {{u_1}^{m+1}_{i+1,j}-2{u_1}^{m+1}_{i,j}+u^{m+1}_{i-1,j}}{\Delta x^2}+\dfrac {{u_1}^{m+1}_{i,j+1}-2{u_1}^{m+1}_{i,j}+{u_1}^{m+1}_{i,j-1}}{\Delta y^2}=f(x_i,y_j), \\[8pt] 
\ \ \ \  i=2,..., I'-1, \ j=2,..., J-1; \ \ {u_1}^{m+1}_{I',j}=N^m_{I',j}, \ j=1, ...,J\\[10pt]
\dfrac {{u_2}^{m+1}_{i+1,j}-2{u_2}^{m+1}_{i,j}+{u_2}^{m+1}_{i-1,j}}{\Delta x^2}+\dfrac {{u_2}^{m+1}_{i,j+1}-2{u_2}^{m+1}_{i,j}+{u_2}^{m+1}_{i,j-1}}{\Delta y^2}=f(x_i,y_j),  \\[8pt]
\ \ \ \  i=I',..., I-1, \ j=2,..., J-1; \ \ {u_1}^{m+1}_{I'-1,j}=N^m_{I'-1,j}, \ j=1, ..., J
\end{array}\right. 
\end {equation} 

\noindent plus the boundary condition ${u_1}^{m+1}_{i,j},\ {u_2}^{m+1}_{i,j}=g(x_i,y_j)$ on $\partial \Omega _1\cap \partial \Omega _2$. The discretized Poisson equation can also be computed via a classic Schwarz iteration algorithm in association with interface conditions: 
\begin {equation} 
\label {classic_interface}
{u_1}_{I',j}^{m+1}={u_2}_{I',j}^m, \ \ {u_2}_{I'-1,j}^{m+1}={u_1}_{I'-1,j}^m, \ \ j=2, ..., J-1
\end {equation} 


For the training by (\ref {interface_solu1}) and (\ref {loss}), the data can be generated using numerical solutions. For instance, let the iterated solutions on the bounaries of the interface zone in the above classic Schwarz iteration and also the value of $g$ on the rest boundaries as the input, and the corresponding solutions on the interfaces to the discretized Poisson equation as the output. Particularly,  
\begin {equation} 
\label{training}
\left\{\begin{array}{llll} 
u_0:      & {u_1}^m_{I'-k-1,j}, \ {u_2}_{I'+k,j}^m,     
& on \ S_1, S_2 \\
          & g(x_i,y_1),\ g(x_i,y_J),  \           
& on \  \partial \Omega _1\cap \partial \Omega _2\\
\bar {s}: & {u_1}^m_{I',j}, {u_2}_{I'-1,j}^m, \         
& on \ \Gamma _1, \Gamma_2 
\end{array}\right . 
\end {equation}

\begin{figure}[] 
  \vspace{-0.5 cm}
  \centering
  \subfloat[][]{\label{Poisson_solu}\includegraphics[width=0.33\textwidth]{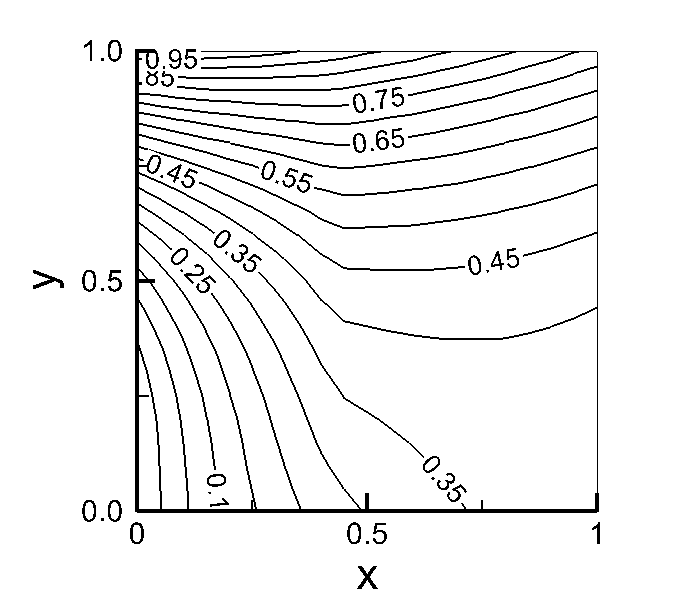}}  \ \ 
  \subfloat[][]{\label{Poisson_error}\includegraphics[width=0.33\textwidth]{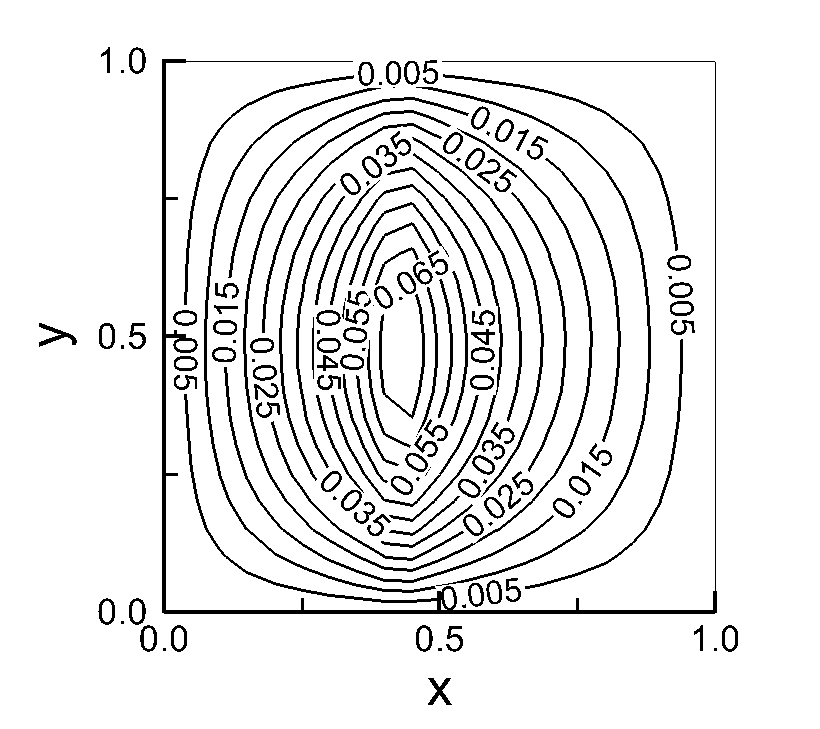}} 
  \vspace{-8pt}
  \caption{\small Computation of coupled Poisson equations.  }
  \label{Poisson_eq_solu}
  \vspace{-8pt}
\end{figure}

Let the grid be $(I,J)=(21,21)$, and interfaces be at $I'-1=20, \ I'=21$, with $k=2$, with a feedorward neural network as the network and Sigamoid function as the activation function \cite {Lagaris1998}. ML Schwarz iterative algorithm (\ref {ML_Schwarz1}) is computed using MatLab, and the solution, $u$ ($=u_1\cup u_2$), is presented in Fig. \ref {Poisson_eq_solu}. For comparison, the solution associated with classic interface condition (\ref {classic_interface}), referred to as a classic numerical solution, $u_{num}$, is also computed. It is seen that the solution obtained with the ML method matches the exact solution. Interestingly, it is noticed that the ML solution may converge faster than the classic numerical solution. For instance, when one hidden layer with 10 neurons is used, the ML method converges in just 18 iterations, while the classic numerical solution requires 121 iterations, with convergence criterion of $10^{-10}$ in iteration residual. 

Now let us examine the error of the ML solution against the exact solution. Since we know that 
\begin {equation}
\label{accuracy}
||u-u_{exact}|| \le ||u-\bar {s}||+||\bar {s}-u_{exact}||
\end {equation} 

\noindent the error can be estimated via two parts. The first part is the error of the ML solution against the training data, and the second part is the error of the training data against the exact solution. The former depends on how the ML solution's structure and training method, i.e., (\ref {interface_solu1}) and (\ref {loss}), and the latter depends on the accuracy order in the discretization of the Poisson equation. These errors are shown in Table \ref {table_Poisson}. Here, and hereafter, in view that an ML run produces results slightly different from those of another, each case's values result from the averages of three runs. The table indicates that the latter is relatively small, and this implies that, in this example, the error of the ML solution primarily comes from its error against the training data. 

\begin{table}[t]
  \centering
  \caption{\small Errors of computation for the Poisson equation. $N_l$ and $N_n$ are respectively the number of layers and neurons. $||\cdot ||$ is infinite norms over $S$, and $Loss$ is in mean square error. } 
  \label{table_Poisson}
  \begin{tabular}{|c|c|c|c|c|}
\hline
       $N_l$, $N_n$  & $||u-u_{exact}||$ & $||u-\bar {s}||$ & $||\bar {s}-u_{exact}||$ & $Loss$ \\ \hline
        1, \ 5  & $6.9765\times 10^{-2}$ &   $6.9797\times 10^{-2}$   & $3.3\times 10^{-5}$      & $4.28\times 10^{-10}$   \\ \hline        
        1, \ 10 & $6.9581\times 10^{-2}$ &   $6.9613\times 10^{-2}$  & $3.3\times 10^{-5}$       & $5.63\times 10^{-10}$   \\ \hline
        1, \ 15 & $6.9524\times 10^{-2}$ &   $6.9557\times 10^{-2}$  & $3.3\times 10^{-5}$       & $4.09\times 10^{-10}$   \\ \hline        
        1, \ 20 & $6.9848\times 10^{-2}$ &   $6.9876\times 10^{-2}$  & $3.3\times 10^{-5}$       & $1.61\times 10^{-10}$   \\ \hline        
        2, \ 10 & $6.9516\times 10^{-2}$ &   $6.9548\times 10^{-2}$  & $3.3\times 10^{-5}$       & $3.14\times 10^{-10}$   \\ \hline
        3, \ 10 & $6.9491\times 10^{-2}$ &   $6.9524\times 10^{-2}$  & $3.3\times 10^{-5}$       & $1.04\times 10^{-9}$   \\ \hline        
  \end{tabular}
\end{table}

\vspace{0.2cm}

\noindent {\bf 3. Advection-diffusion equations} 

\vspace{0.1cm}

Consider an initial value problem of two coupled advection-diffusion equations over subdomains  $\Omega _1$ and $\Omega _2$, which overlap with each other with interfaces $\Gamma _1$ and $\Gamma _2$, respectively. When marching from time level $n$ to $n+1$ ($n=0,1,2, ..., N$) within time interval $[0,T]$, the coupled problem is solved by an ML Schwarz waveform relaxation algorithm:

\begin{equation}
\label{ML_Schwarz2}
\begin{array}{ll}
\left\{\begin{array}{ll}
{u_1}_t^{n+1}+a_1\nabla \cdot {u_1}^{n+1}=b_1\Delta {u_1}^{n+1}, & \ t, \textbf{x}\in (0,T]\times \Omega _1 \\
{u_1}^{n+1}=s^n, & \ t, \textbf{x} \in (0,T]\times \Gamma _1 \\
{u_1}^0=g, & \ t,\textbf{x}\in (t=0)\times \Omega _1
\end{array}\right . \\[20pt]
\left\{\begin{array}{ll}
{u_2}_t^{n+1}+a_2\nabla \cdot {u_2}^{n+1}=b_1\Delta {u_1}^{n+1}, & \ t, \textbf{x}\in (0,T]\times \Omega _2 \\
{u_2}^{n+1}=s^n, & \ t, \textbf{x} \in (0,T]\times \Gamma _2 \\
{u_2}^0=g, & \ t,\textbf{x}\in (t=0)\times \Omega _2
\end{array}\right .
\end{array}
\end{equation}

\noindent where $g \in C^0$, and $s$ is an ML solution that provides transmission conditions. As an important feature in this algorithm, Schwarz iteration between the two subdoamins is not necessary when marching from time level $n$ to $n+1$, since the ML may predict the solution on the interfaces at time level $n+1$. But, if needed, Schwarz iteration may also be made, with $n$ being replaced by $m$, the index of iteration, in the above algorithm. 

Similarly to the situation of the above Poisson equation, let us introduce an interface zone, $S$, that encloses the interfaces (e.g., Fig. \ref {domain_division}). Let $s$ be the ML solution, $u$ be the solution of the corresponding classic initial value problem, and $\textbf {u}_0 $ be prescribed initial condition. Then, the problem of the ML solution becomes an initial and boundary problem: Find $\mathcal {M}$ such that

\begin {equation} 
\label{ML_problem2}
\left\{\begin{array}{ll} 
s= \mathcal {M}(u_0)\approx u,  \ \ \ & t,\textbf{x}\in (t_0,t]\times S \\
s=u_0, \ \ \ & t,\textbf{x}\in (t=t_0)\times S, S_1\cup S_2\times (t_0,t], \partial \Omega _1\cap \partial \Omega _2\times (t_0,t]
\end{array}\right . 
\end {equation}

\noindent where $t_0$ is a prescribed value of time. Now $\mathcal {M}$ may be considered as an operator that mimics the advection-diffusion equations. 

\begin{wrapfigure}{r}{0.4\textwidth} 
\vspace {-0.5 cm}
\centering
\includegraphics[width=0.99\linewidth]{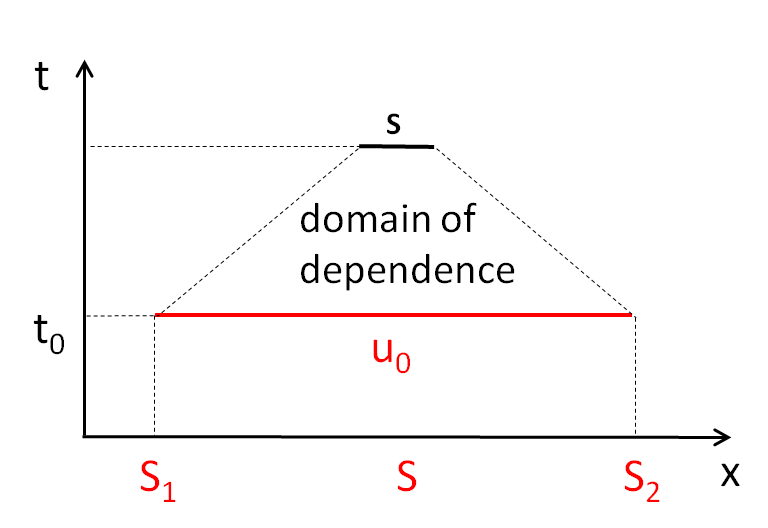} 
\vspace {-0.5 cm}
\caption{\small Schematic representation of the local approach for the training data.}
\label{ANN_training}
\vspace {-0. cm}
\end{wrapfigure}
In ML problem (\ref {ML_problem2}), $u_0$ is located in the interface zone and at its boundaries. Three training approaches in terms of the place where $u_0$ resides are possible: i) global: $ \Omega \times (t=0)$, 
ii) regional: $S\times (0, T]$, $S_1\cup S_2\times (0,T] $, and iii) local: $S\times (t=t_0) $. In the global approach, the input data, $u_0$, is actually the initial value of the physical problem itself, and the ML solution needs no information from numerical solutions in the subdomains, actually decoupling the numerical solutions completely. The local approach corresponds to a local initial value problem over the interface zone. In the situation shown in Fig. \ref {ANN_training}, the domain of dependence of $s$ is the trapezoid with the interface zone as its base (in red), where $u_0$ is located. It is anticipated that, for the seek of stability, the size of the interface zone should be wide enough so this domain contains the domain of dependence of the solution to problem (\ref {ML_Schwarz2}). But, the size is not necessarily very big since the solution far away has little effects. When the local training approach in adopted, the numerical solutions in the two subdomains are coupled in two way with the ML solution.  

Now let us do numerical experiments on the ML-based coupling. 
First, consider two identical one-dimensional equations in two subdomains. Particularly, $a_1,a_2=a=1, b_1,b_2=b=0.1; \ g=-sin(\pi x), \ -1<x<1; \ u=0, \ x=-1, \ 1$. The problem has an analytical solution \cite {Abdelkader2015}: 
\begin {equation} 
\label{analytical_solu}
\begin{array}{ll}
u =& 16\pi^2 ab^3 e^{(x-ct/2)a/2b } \\
& \times \left[ 
sinh(\dfrac{a}{2b}) \sum_{p=0}^{\infty}\dfrac{(-1)^p\ 2p\ sin(p\pi x)\ e^{-b p^2 \pi^2 t}}
{a^4+8(ab\pi )^2(p^2+1)+16(\pi b)^4(p^2-1)^2} \right.\\ 
& \quad \left.+
cosh(\dfrac{a}{2b}) \sum_{p=0}^{\infty}\dfrac{(-1)^p\ (2p+1)\ cos((2p+1)/(2\pi x))\ e^{-(2p+1)^2b\pi^2 t/4}}
{a^4+(ab\pi)^2(8p^2+8p+10)+(\pi b)^4(4p^2+4p-3)^2} \right]
\end{array} 
\end {equation} 

On grid $x_1,x_2,...,x_I$ ($x_1=-1, \ x_2=1$), let $\Omega _1$ be $\ x_1<x<x_I'$, and $\Omega _2$ be $x_{I'-1}<x<x_I$, with their interfaces at $x_{I'-1},x_{I'}$. Then problem (\ref {ML_Schwarz2}) is discretized

\begin {equation} 
\left\{\begin{array}{ll}
\dfrac {{u_1}^{n+1}_i-{u_1}^n_i}{\Delta t}
+a_1\dfrac {{u_1}^{n+1}_{i+1}-{u_1}^{n+1}_{i-1}}{2\Delta x}
=b_1\dfrac {{u_1}^{n+1}_{i+1}+{u_1}^{n+1}_{i+1}-2{u_1}^{n+1}_i}{\Delta x^2}, \\[8pt]
\ \ \ \  i=2,..., I'-1, \ j=2,..., J-1; \ \ {u_1}_{I'}^{n+1}=N^n_{I'} \\[10pt]
\dfrac {{u_2}^{n+1}_i-{u_2}^n_i}{\Delta t}
+a_2\dfrac {{u_2}^{n+1}_{i+1}-{u_2}^{n+1}_{i-1}}{2\Delta x}
=b_2\dfrac {{u_2}^{n+1}_{i+1}+{u_2}^{n+1}_{i+1}-2{u_2}^{n+1}_i}{\Delta x^2}, \\[8pt]
\ \ \ \  i=I',..., I-1, \ j=2,..., J-1; \ \ {u_2}_{I'-1}^{n+1}=N^n_{I'-1}
\end{array}\right. 
\end {equation} 

\noindent A corresponding classic Schwarz iteration algorithm uses the following interface conditions: 

\begin {equation} 
{u_1}_{I'}^{n+1}={u_2}_2^n, \ \ {u_2}_1^{n+1}={u_1}_{I'-1}^n
\end {equation} 

Let us train the network via (\ref {interface_solu1}), (\ref {loss}), and the global approach; the initial and boundary values of the problem are the input, and above analytical solution is the output of the network. Particularly, 

\begin {equation} 
\label{global_training}
\left\{\begin{array}{llll} 
&u_0: & sin(\pi x_i),                  & on \ (t=0)\times \Omega , \\ 
&     & -1, \mbox {or}, 1,             & on \ (0,1]\times \partial \Omega \\
&\bar {s}: &u(t^l,x_{I'-1}), u(t^l,x_{I'}), & on \ (0,1]\times \Gamma _1, \ (0,1]\times \Gamma _2 
\end{array}\right . 
\end {equation}

\noindent in which $l=1,2,...L$, with step $\Delta t'$. The ML solution is constructed as a function of $N$: 

\begin {equation} 
\label{interface_solu2}
\left\{\begin{array}{llll} 
s=A(t,x)+ t(1-t)(x+1)(x-1)N(t,\textbf {x}, u_0, \boldsymbol {\alpha }) \\
A(t,x)= (1-t)\left (u(0,x)-\dfrac {1-x}{2}u(0,-1)-\dfrac {1+x}{2}u(t,0)\right )\\
\ \ \ \ \ \ \ \ \ \  +t \left (u(1,x)-\dfrac {1-x}{2}u(1,-1)-\dfrac {1+x}{2}u(1,1)\right ) \\
\ \ \ \ \ \ \ \ \ \  +\dfrac {1-x}{2}u(-1,t)+\dfrac {1+x}{2}u(1,t)
\end{array}\right . 
\end {equation}

\noindent which automatically satisfies the initial and boundary conditions. This reduces the original constrained optimization problem to an unconstrained one in training the network \cite {Lagaris1998}. 

\begin{figure}[t] 
  \vspace{-0.5 cm}
  \centering
  \subfloat[][]{\label{identical_pde_solu}\includegraphics[width=0.33\textwidth]{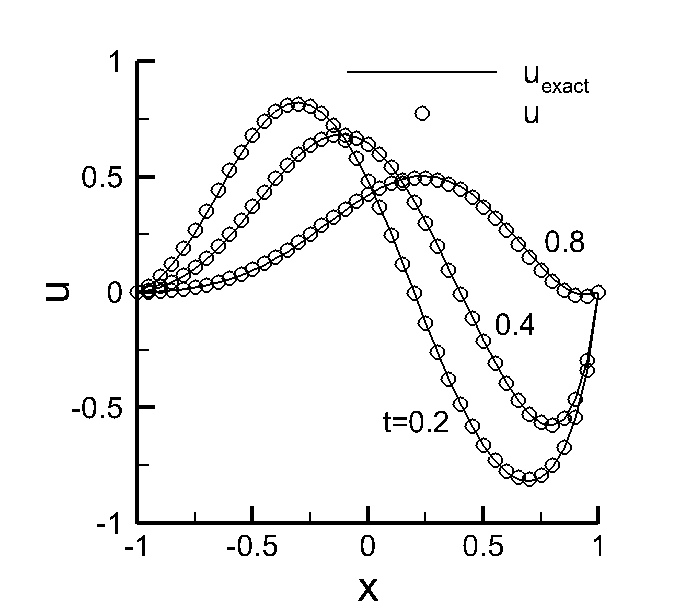}}  \ \ 
  \subfloat[][]{\label{identical_pde_error}\includegraphics[width=0.33\textwidth]{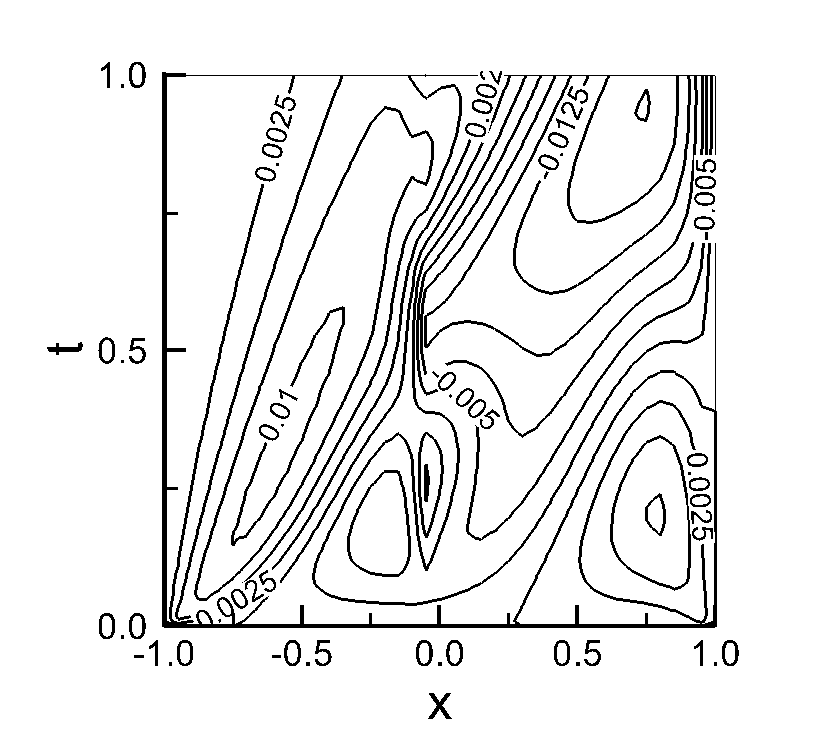}} 
  \vspace{-8pt}
  \caption{\small computation of advection-diffusion equations with identical coefficients. 1) Computed and exact solutions at different time. 2) Error distribution of $u-u_{exact}$.  }
  \label{advection_diffusion_solu}
  \vspace{-8pt}
\end{figure}

In computation, let $(I,N,L)=(41,101,32)$, and use the network and activation function as above. The computed and exact solutions, plus the error distribution, are shown in Fig. \ref {advection_diffusion_solu}.  The errors are presented in Table \ref {table_advection_diffusion_same_coef}. From the table, it is seen that, although not substantially, the solution exhibits convergence with regard to numbers of neurons when one hidden layer is used, which is consistent with an earlier analysis \cite {Sirignan2018}, while the training loss decreases dramatically. Note that the training mesh is not the same as the mesh for the numerical solution, and thus the ML solution does not simply repeat the analytical solution at the interfaces in the figure. 

\begin{table}[h!]
  \centering
  \caption{\small Errors of computation for the advection-diffusion equations with identical coefficients. $N_l$ and $N_n$ are respectively the number of layers and neurons. $||\cdot ||$ is infinite norms over $\Omega \times [0,T]$, and $Loss$ is in mean square error.} 
  \label{table_advection_diffusion_same_coef}
  \begin{tabular}{|c|c|c|c|c|}
\hline
       $N_l$, $N_n$     & $||u-u_{exact}||$           & $Loss$                   \\ \hline
        1, \  5           & $2.35\times 10^{-2}$      &  $7.37\times 10^{-1}$     \\ \hline       
        1, \  10          & $1.70\times 10^{-2}$      &  $9.04\times 10^{-2}$     \\ \hline
        1, \  15          & $1.65\times 10^{-2}$      &  $1.28\times 10^{-2}$     \\ \hline
        1, \  20          & $1.64\times 10^{-2}$      &  $3.36\times 10^{-3}$     \\ \hline
        2, \  10          & $1.67\times 10^{-2}$      &  $7.90\times 10^{-7}$     \\ \hline
        3, \  10          & $1.67\times 10^{-2}$      &  $9.04\times 10^{-9}$     \\ \hline
  \end{tabular}
\end{table}

\begin{figure}[t] 
  \vspace{-0.5 cm}
  \centering
  \subfloat[][]{\label{different_pde_solu}\includegraphics[width=0.33\textwidth]{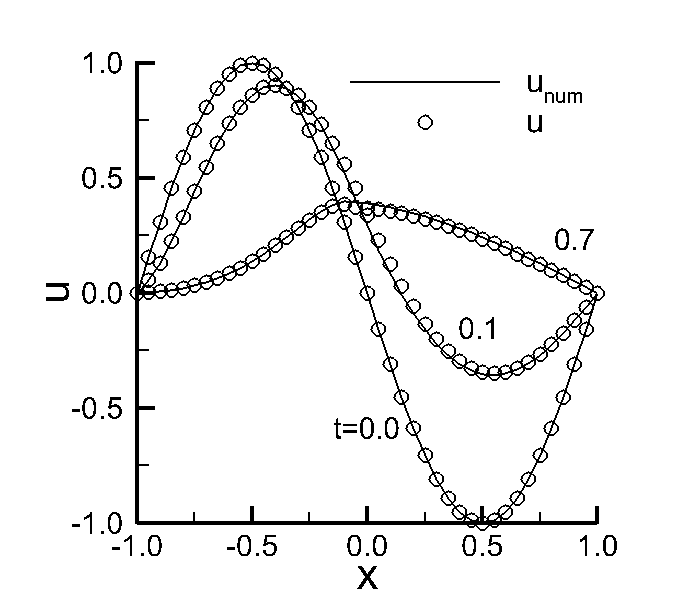}}  \ \ 
  \subfloat[][]{\label{different_pde_error}\includegraphics[width=0.33\textwidth]{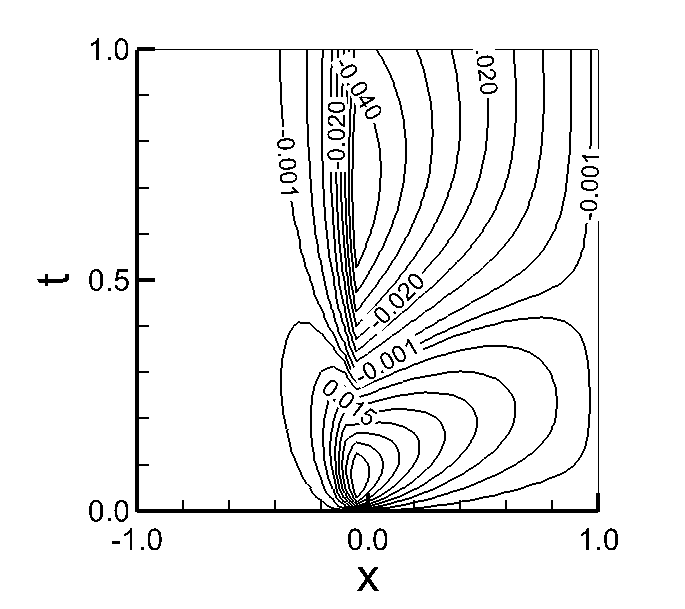}} 
  \vspace{-8pt}
  \caption{\small computation of advection-diffusion equations with different coefficients. 1) ML Computed and classic numerical solutions at different time. 2) Error distribution of $u-u_{num}$.  }
  \label{different_pde}
  \vspace{-8pt}
\end{figure}

Second, consider a heterogeneous problem, that is, a situation of equations with distinct coefficients in the two subdomains: $a_1=1, \ a_2=0.1, \ b_1=0.1,\ b_2=1$.  In this situation, no analytical solution is available, and the numerical solution obtained with the classic Schwarz iteration is used to train the network. The local approach for the training is adopted, with the input and output as follows:  
\begin {equation} 
\label{local_training}
\left\{\begin{array}{llll} 
u_0: & {u_1}^n_{I'-k-1}, ..., {u_1}^n_{I'-1},{u_2}^n_{I'},...,{u_2}^n_{I'+k}, \ & on \ (t=0)\times S , \\
\bar {s}: &u^{n+1}_{I'-1}, u^{n+1}_{I'}, \ &  on \ (t=\Delta t)\times \Gamma _1, \ (t=\Delta t)\times \Gamma _2
\end{array}\right . 
\end {equation}

\noindent for all $n=1,2, ....N$. Note that the mesh for the training and that for the numerical solution associated with the ML solution are the same as those above. The computed solutions are presented in Fig. \ref {different_pde} and Table \ref {table_pde_diff_coef}, from which it is seen that the ML numerical solution behaves similarly to those with same coefficients as discussed above.  


\begin{table}[h!]
  \centering
  \caption{\small Errors of computation for the advection-diffusion equations with different coefficients. 
$||\cdot ||$ is infinite norms over $S \times [0,T]$, and $Loss$ is in mean square error.} 
  \label{table_pde_diff_coef}
  \begin{tabular}{|c|c|c|c|c|}
\hline
       $N_l$, $N_n$     & $||u-u_{num}||$           & $Loss$                   \\ \hline
        1, \  2           & $6.90\times 10^{-2}$      &  $1.01\times 10^{-10}$     \\ \hline       
        1, \  4           & $5.02\times 10^{-2}$      &  $6.95\times 10^{-11}$     \\ \hline
        1, \  8           & $4.60\times 10^{-2}$      &  $3.33\times 10^{-11}$     \\ \hline
        1, \  16          & $4.35\times 10^{-2}$      &  $3.02\times 10^{-11}$     \\ \hline
        2, \  4           & $4.88\times 10^{-2}$      &  $1.90\times 10^{-11}$     \\ \hline
        3, \  4           & $4.50\times 10^{-2}$      &  $1.71\times 10^{-11}$     \\ \hline
  \end{tabular}
\end{table}

\vspace{0.2cm}

\noindent {\bf 4. Discussion}

\vspace{0.1cm}

In ML coupling frameworks (\ref {ML_Schwarz1}) and (\ref {ML_Schwarz2}), numerical solutions in subdomains are coupled via an ML solution in between. This is different from conventional coupling methods, in which the numerical solutions are coupled directly with each other. It is expected that such a difference could enable the ML-based methods to exhibit advantages over the conventional methods in the following aspects: 

{\it Capability} --- An ML solution within the interface zone introduces a `buffer zone' between solutions of PDEs in subdoamins, and thus may relax the rigidity of well-posedness for their coupling, which could be hard to achieve but is necessary for a conventional coupling approach. As long as the ML solution assures that the PDE problems in individual sudomains are well-posed, which we know more and is usually clear, the computation of the coupled PDEs tends to be good. As a result, an ML approach may avoid difficulties/singularities encountered in a conventional approach, e.g., non-physical solutions at interfaces. 

{\it Efficiency} --- An ML interface algorithm may be trained offline, and also it does not need Schwarz iteration, which requires multiple times of computation of PDEs in each subdomain and is time-consuming. These lead to the potential to reduce the computation time of coupled PDEs significantly. 
 
{\it Universality} --- It is conjectured that an ML-based coupler can be `universal': a coupler developed for one set of PDE solvers is applicable to coupling another set of solvers with the same PDEs but different discretization and software. Traditional coupling methods heavily depend on discretization and software and lack such `universality'.

The approach of coupling PDEs via ML is built on the capabilities of ML methods, which are anticipated to bear solid foundation in mathematics. ML is not merely interpolation or data fitting, and its performance in prediction has been proven promising, especially in the case of deep learning built on PDEs \cite {Lee2017,Wiewel2019}. Such performance is being recognized and agreed upon in the communities. For this reason, research has been initiated on the mathematical foundation and mechanism of ML \cite {MoDL2020, DOE2020}. Moreover, our preliminary results reported in this paper have shown that ML coupling works as intended.  

ML coupling of PDEs is basically an untapped topic, and it needs further exploration. This study presents the ideas, concepts, and methods for ML-based coupling of PDEs, and, as first-hand knowledge, it provides numerical examples of its feasibility and performance. This study is exploratory and preliminary, and further research is necessary on the ML coupling to explore the above listed potential advantages via systematic numerical and theoretical analyses. Nevertheless, this short communication paper indicates that the ML coupling is promising, and, hopefully, it will attract the community's attention to the topic and its further study. 

\vspace{0.2cm}

\noindent {\bf Acknowledgments.} HT, YJ, and WD are supported by NSF (DMS \#1622453, \#1622459).

\bibliography{bbi_ML_ddm}

\begin{thebibliography}{10}

\bibitem{Schwarz1869}
H.~A. Schwarz.
\newblock \"{U}ber einige abbildungsaufgaben.
\newblock {\em Ges. Abh.}, 11:65--83, 1869.

\bibitem{Califano2018}
G.~Califano and D.~Conte.
\newblock Optimal {Schwarz} waveform relaxation for fractional diffusion-wave
  equations.
\newblock {\em Applied Numerical Mathematics}, 127:125–141, 2018.

\bibitem{Gastaldi1989}
F.~Gastaldi and A.~Guarteroni.
\newblock On the coupling of hyperbolic and parabolic systems: Analytical and
  numerical approach.
\newblock {\em Applied Numerical Mathematics}, 6:3--31, 1989/90.

\bibitem{Gander1998}
M.~J. Gander and A.~M. Stuart.
\newblock Space-time continuous analysis of waveform relaxation for the heat
  equation.
\newblock {\em SIAM J. Sci. Comput.}, 19:2014--2031, 1998.

\bibitem{Martin2004}
V.~Martin.
\newblock An optimized schwarz waveform relaxation method for the unsteady
  convection diffusion equation in two dimensions.
\newblock {\em Computers \& Fluids}, 33:829–837, 2004.

\bibitem{Gander2007}
M.~J. Gander and L.~Halpern.
\newblock Optimized {Schwarz} waveform relaxation methods for advection
  diffusion problems.
\newblock {\em SIAM J. Numer. Anal.}, 45:666--697, 2007.

\bibitem{Berger1987}
M.~J. Berger.
\newblock On conservation at grid interfaces.
\newblock {\em SIAM J. Numer. Anal.}, 24:967--984, 1987.

\bibitem{Tang1999}
H.~S. Tang and T.~Zhou.
\newblock On non-conservative algorithms for grid interfaces.
\newblock {\em SIAM J. Numer. Anal.}, 37:173--193, 1999.

\bibitem{keyes2012}
D.~E. Keyes.
\newblock Multiphysics simulations: Challenges and opportunities.
\newblock {\em Int. J. High Performance Computing Application}, 27(1):4--86,
  2012.

\bibitem{Tang2020}
H.S. Tang, R.D. Haynes, and G.~Houzeaux.
\newblock A review of domain decomposition methods for simulation of fluid
  flows: Concepts, algorithms, and applications.
\newblock {\em Arch Computat Methods Eng}, 2020.

\bibitem{Minami:2012jp}
S.~Minami, H.~Kawai, and S.~Yoshimura.
\newblock A monolithic approach based on the balancing domain decomposition
  method for acoustic fluid-structure interaction.
\newblock {\em Journal of Applied Mechanics}, 79(1):010906, 2012.

\bibitem{Tang2014a}
H.~S. Tang, K.~Qu, and X.~G. Wu.
\newblock An overset grid method for integration of fully {3D} fluid dynamics
  and geophysical fluid dynamics models to simulate multiphysics coastal ocean
  flows.
\newblock {\em J. Comput. Phys.}, 273:548–571, 2014.

\bibitem{blayo2016}
E.~Blayo.
\newblock About interface conditions for coupling hydrostatic and
  nonhydrostatic navier-stokes flows.
\newblock {\em Discrete and Continuous Dynamical Systems Series, S},
  9(5):1565--1574, 2016.

\bibitem{Qu2019a}
K.~Qu, H.S. Tang, and A.~Agrawal.
\newblock Integration of fully {3D} fluid dynamics and geophysical fluid
  dynamics models for multiphysics coastal ocean flows: Simulation of local
  complex free-surface phenomena.
\newblock {\em Ocean Modelling}, 135:14--30, 2019.

\bibitem{Eisenmann2018}
M.~Eisenmann and E.~Hansen.
\newblock Convergence analysis of domain decomposition based time integrators
  for degenerate parabolic equations.
\newblock {\em Numerische Mathematik}, 140:913–938, 2018.

\bibitem{Part1994}
E.~P{\"a}rt-Enander and B.~Sj{\"o}green.
\newblock , conservative and nonconservative interpolation between overlapping
  grids for nite volume solutions of hyperbolic problems.
\newblock {\em Comput. \& Fluids}, 23:551--574, 1994.

\bibitem{Wu1996}
Z.~L. Wu.
\newblock On uniqueness of steady state solutions for difference equations on
  overlapping grids.
\newblock {\em SIAM J. of Numer. Anal.}, 33(4):1336--1357, 1996.

\bibitem{Tang1996}
H.~S. Tang, D.~L. Zhang, and C.~H. Lee.
\newblock Comments on algorithms for grid interfaces in simulating {Euler}
  flows.
\newblock {\em Comm. Nonlinear Science and Numerical Simulation}, 1:50--54,
  1996.

\bibitem{Tang2016a}
H.~S. Tang, K.~Qu, X.~G. Wu, and Z.~K. Zhang.
\newblock {\em Domain decomposition for a hybrid fully {3D} fluid dynamics and
  geophysical fluid dynamics modeling system: A numerical experiment on a
  transient sill flow.}
\newblock Domain Decomposition Methods in Science and Engineering XXII,
  407-414. Lecture Notes in Computational Science and Engineering, XXII. Ed:
  Dickopf, T and Gander, MJ and Halpern, L and Krause, R and Pavarino, LF.
  Springer, 2016.

\bibitem{Tang2019}
H.~S. Tang and Y.~J. Liu.
\newblock {\em Coupling of Navier-Stokes Equations and Their Hydrostatic
  Versions for Ocean Flows: A Discussion on Algorithm and Implementation. Proc.
  DD25.}
\newblock Domain Decomposition Methods in Science and Engineering, Lecture
  Notes in Computational Science and Engineering. Springer, In press.

\bibitem{Tang2003}
H.~S. Tang, C.~Jones, and F.~Sotiropoulos.
\newblock An overset grid method for 3d unsteady incompressible flows.
\newblock {\em J. Comput. Phys.}, 191:567--600, 2003.

\bibitem{Chen2014b}
R.~L. Chen, Y.~Q. Wu, Z.Z. Yan, Y.~Zhao, and X.-C. Cai.
\newblock A parallel domain decomposition method for 3d unsteady incompressible
  flows at high {Reynolds} number.
\newblock {\em J Sci Comput}, 58:275--289, 2014.

\bibitem{Yassin2018}
N.~I.~R. Yassin, S.~Omran, E.~M.~F. El~Houby, Enas M., and H.~Allam.
\newblock Machine learning techniques for breast cancer computer aided
  diagnosis using different image modalities: A systematic review.
\newblock {\em Computer Methods and Programs in Biomedicine}, 156:25--45, 2018.

\bibitem{Park2018}
W.~Park, S.~Kim, K.~L. Kim, and J.~Kim.
\newblock Alphago's decision making.
\newblock {\em J. Logics and Their Applications}, 6:105--155, 2019.

\bibitem{Tompson2017}
J.~Tompson, K.~Schlachter, P.~Sprechmann, and K.~Perlin.
\newblock Accelerating eulerian fluid simulationwith convolutional networks.
\newblock {\em PMLR}, 70:3424--3433, 2017.

\bibitem{Raissi2019}
M.~Raissi, Z.~Wang, and G.~E. Triantafyllou, M. S.and~Karniadakis.
\newblock Deep learning of vortex-induced vibrations.
\newblock {\em J. Fluid Mech.}, 861:119--137, 2019.

\bibitem{Lee2017}
S.~Lee and D.~You.
\newblock Prediction of laminar vortex shedding over a cylinder using deep
  learning. arxiv:1712.07854 [physics.flu-dyn].
\newblock {\em J. Fluid Mechanics, submitted}, 2017.

\bibitem{Wiewel2019}
S.~Wiewel, M.~Becher, and N.~Thuerey.
\newblock Latent space physics: Towards learning the temporal evolution of
  fluid flow.
\newblock {\em Eurographics}, 38, 2019.

\bibitem{Sirignan2018}
J.~Sirignan and K.~Spiliopoulos.
\newblock {DGM}: A deep learning algorithm for solving partial differential
  equations.
\newblock {\em J. Comput. Phys.}, 375:1339--1364, 2018.

\bibitem{matlabweb}
MathWorks.
\newblock Matlab.
\newblock https://www.mathworks.com/products/matlab.html.

\bibitem{pytorchweb}
{PyTorch}.
\newblock From research to production.
\newblock https://pytorch.org/.

\bibitem{tensorflowweb}
{Tensorflow}.
\newblock An end-to-end open source machine learning platform.
\newblock https://www.tensorflow.org/.

\bibitem{Pawar2020}
S.~Pawar, S.~E. Ahmed, and M.~San.
\newblock Interface learning in fluid dynamics: statistical inference of
  closures within micro-macro coupling models. arxiv:2008.04490v1
  [physics.comp-ph] 11 aug 2020, 2020.

\bibitem{Tang2019b}
H.~S. Tang, Y.~J. Liu, and M.~Grossberg.
\newblock Research proposals. 1) collaborative research: {CDS\& E}: {A}
  modeling system and simulation for emerging complex flow problems in coastal
  waters. {NSF, Sept., 2019}; 2) {A} machine-learning framework to couple
  numerical solutions of partial differential equations towards flow
  simulations, {DOE, Apr. 2020}.

\bibitem{Toselli2005}
A.~Toselli and O~Widlund.
\newblock {\em Domain Decomposition Methods – Algorithms and Theory}.
\newblock Springer, 2005.

\bibitem{Lagaris1998}
I.~E. Lagaris, A.~Likas, and D.~I. Fotiadis.
\newblock Artificial neural networks for solving ordinary and partial
  differential equations.
\newblock {\em {IEEE} Transactions on Neural Networks}, 9:987--1000, 1998.

\bibitem{Baymani2015}
M.~Baymani, S.~Effati, H.~Niazmand, and A.~Kerayechian.
\newblock Artificial neural network method for solving the {Navier–Stokes}
  equations.
\newblock {\em Neural Comput \& Applic}, 26:765–773, 2015.

\bibitem{Kou2019}
J.~Kou and W.~Zhang.
\newblock Multi-fidelity modeling framework for nonlinear unsteady aerodynamics
  of airfoils.
\newblock {\em Applied Mathematical Modelling}, 76:832–855, 2019.

\bibitem{Liu2019}
D.~Liu and Y.~Wang.
\newblock Multi-fidelity physics-constrained neural network and its application
  in materials modeling.
\newblock {\em J. Mech. Design}, 141(12):121403, 2019.

\bibitem{Abdelkader2015}
Abdelkader Mojtabi and Michel Deville.
\newblock One-dimensional linear advection–diffusion equation:analytical and
  finite element solutions.
\newblock {\em Computers \& Fluids}, 107:189--195, 2015.

\bibitem{MoDL2020}
{NSF}.
\newblock {NSF-Simons} research collaborations on the mathematical and
  scientific foundations of deep learning ({MoDL}), {Division of Mathematical
  Sciences}, solicitation 20-540, 2020.

\bibitem{DOE2020}
{DOE, Advanced Scientific Computing Research}.
\newblock Scientific machine learning for modeling and simulations:
  {DE-FOA-2319}, 2020.

\end{thebibliography}
\end{document}